\documentclass[12pt]{article}
\usepackage{latexsym,amsmath,theorem,amssymb,biometri,graphicx,mathrsfs,
  setspace}

\let\bibhang\relax
\usepackage{natbib}

\usepackage{url}

\usepackage{rotating}

\allowdisplaybreaks

\newcommand{\nc}{\newcommand}

\newcommand*\xbar[1]{%
  \hbox{%
    \vbox{%
      \hrule height 0.5pt
      \kern0.275ex
      \hbox{%
        \kern0em
        \ensuremath{#1}%
        \kern0em
      }%
    }%
  }%
}

\nc{\w}{\wedge}
\nc{\h}{\widehat}
\nc{\td}{\widetilde}
\nc{\T}{\top}

\nc{\e}{E}
\nc{\he}{\h{\e}}

\nc{\al}{\pmb{\eta}}
\nc{\ha}{\h{\al}}
\nc{\ba}{\xbar{\al}}
\nc{\ta}{\td{\al}}

\nc{\bt}{\pmb{\beta}}
\nc{\hb}{\h{\bt}}
\nc{\tb}{\td{\bt}}

\nc{\epi}{\varepsilon}
\nc{\dconv}{\stackrel{d}{\rightarrow}}

\DeclareFontFamily{U}{mathx}{}
\DeclareFontShape{U}{mathx}{m}{n}{<-> mathx10}{}
\DeclareSymbolFont{mathx}{U}{mathx}{m}{n}
\DeclareMathAccent{\widecheck}{0}{mathx}{"71}

\nc{\cb}{\widecheck{\bt}}

\theoremstyle{bkaasp}
\newtheorem{condition}{Condition}
\theoremstyle{bkathm}
\newtheorem{theorem}{Theorem}

\newtheorem{coroll}[theorem]{Corollary}
\newtheorem{prop}[theorem]{Proposition}

\newtheorem{remark}{Remark}

\begin{document}
\baselineskip=25pt

\begin{center}
{\Large\bf Tightening Control in Neyman--Pearson Linear Classification}\\
      {\large\sc Yijian Huang}\\
      Department of Biostatistics and Bioinformatics,
  Emory University,\\
  Atlanta, Georgia 30322, U.S.A.\\
  yhuang5@emory.edu
\end{center}
\vspace*{.1in}

\section*{Abstract}
  Neyman--Pearson classification prioritizes one class
  by constraining its accuracy above a prespecified level, and then takes
  the accuracy of the other class as the utility objective.
  This paradigm is well suited for disease screening and diagnosis, among
  other applications. Statistical learning under this framework is
  complicated since classifier performance determines
  its acceptability. Furthermore, no learned classifier that is consistent
  for the oracle classifier can guarantee satisfaction of the control
  constraint in finite samples. Classical learning theory targets a
  control-relaxed empirical utility maximization (EUM) classifier. However,
  even the EUM classifier fails to achieve the desired control level on
  average. We conjecture that this under-control phenomenon is a manifestation
  of the over-optimism bias well known in standard statistical learning,
  and develop asymptotic theory to confirm it. Motivated by this insight, we
  propose refined learning procedures under two accuracy control
  strategies for the prioritized class: one controlling accuracy in
  expectation and the other with high probability. We further develop
  training-data-based methods to predict and infer class-specific accuracies
  of the resulting classifiers. Simulation studies demonstrate favorable
  finite-sample performance, and we illustrate the proposed methods with an
  application to cancer detection.

\noindent KEY WORDS:
  Cross-audit projection; Cube root asymptotics; Nonparametric classification;
  Over-optimism; Second-order asymptotics; Statistical learning.

\vspace*{.3in}

\section{Introduction}

In binary classification, misclassification costs are often not only
asymmetric across classes but also difficult to quantify. A prominent
example is disease screening and diagnosis, a critical component of precision
medicine: false negatives and false positives typically carry different
clinical consequences, some of which---such as anxiety and pain---resist
precise cost quantification. Common performance metrics such as
overall misclassification rate and expected misclassification cost are
therefore either inappropriate or infeasible in this context. Clinical
practice instead prioritizes one class by constraining its accuracy above a
specified level, and then takes the accuracy of the other class as the
utility objective \citep[e.g.,][]{catalona98,sanda,skates}. The
corresponding performance metric is thus specificity at a controlled
sensitivity, or sensitivity at a controlled specificity, depending on which
class is prioritized \citep[e.g.,][]{huanga}. This type of
classification---analogous to maximizing power at a fixed type~I error level
in hypothesis testing---has a long history in the statistics literature
\citep[e.g.,][]{greenhouse}, and has more recently been termed the
Neyman--Pearson paradigm in the machine learning literature
\citep[e.g.,][]{tong16}. As clinical tests often leverage multiple features
for improved performance, this article focuses on the fundamental problem of
linear classification under this paradigm.

One learning strategy separates the estimation of the linear feature
combination from threshold estimation. For feature combination, parametric
methods such as linear discriminant analysis and logistic regression are
commonly used but may be suboptimal under model misspecification. Popular
nonparametric methods, such as empirical maximization of the area under the
receiver operating characteristic curve \citep{pepe06}, can also be suboptimal
because they are not designed to optimize the performance metric of interest.
Once the feature combination has been estimated, threshold estimation
addresses the Neyman--Pearson criterion. If an independent dataset is reserved
for threshold estimation, methods developed for a single feature
\citep[e.g.,][]{greenhouse,tong18,huanga} can be applied by treating the
estimated feature combination as fixed. However, this approach is
statistically inefficient due to the additional data requirement. Without a separate
dataset, existing alternatives
are limited to specific parametric models \citep[e.g.,][]{wang}.

\begin{figure}[t]
  \centerline{\includegraphics[width=5.5in]{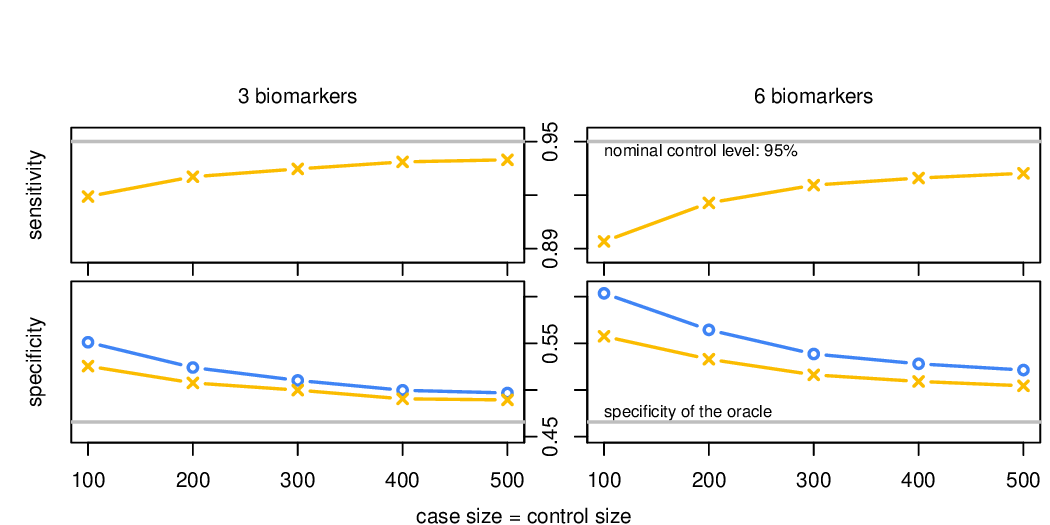}}
  {\caption{\label{fig1}Simulation results for the EUM classifier obtained
      under Scenario~A of Section~5.1, targeting maximum specificity at
      controlled 95\% sensitivity in the cases of 3 and 6 biomarkers. The top
      panels show averaged prediction sensitivity ($\times$), while empirical
      sensitivity equals the nominal control level by construction. The
      bottom panels display averaged prediction ($\times$) and
      empirical ($\circ$) specificities.}}
\end{figure}

A more direct approach to estimating the oracle Neyman--Pearson linear
classifier is to maximize the empirical utility or one of its variants.
General learning theory has been developed for Vapnik--Chervonenkis (VC)
classes of classifiers, with linear classifiers as a special case; see
\cite{cannon} and \cite{scott05}, among others. To obtain nonasymptotic
guarantees on excess risk, however, these methods target a relaxed control
level. Moreover, their theoretical guarantees have limited practical
relevance, and the associated computational challenges remain unresolved.
\cite{meisner} proposed maximizing a kernel-smoothed utility estimate, but
the statistical properties of the resulting classifier remain largely
unexplored beyond consistency. More recently, \cite{huangb} developed a
practical algorithm for empirical utility maximization (EUM) and established
cube-root asymptotics \citep{kim} for the estimated combination coefficient.
Nevertheless, in finite samples, the EUM classifier systematically fails to
attain the nominal control level for the prioritized class.
Figure~\ref{fig1} illustrates this phenomenon using a simulation study
described in Section~5.1: When specificity is optimized subject to a nominal
sensitivity of 95\%, the prediction sensitivity falls below the target level
on average.

Neyman--Pearson classification has a distinctive feature: no estimated
classifier that is consistent for the oracle can guarantee satisfaction of
the control constraint in finite samples. Two main strategies have been
proposed to address this issue. \cite{greenhouse}, \cite{huanga},
\cite{meisner}, and \cite{huangb} implicitly target a control level that is
unbiased for the nominal value. In contrast, \cite{rigollet}, \cite{tong18},
\cite{tong20}, and \cite{wang} focus on ensuring that the nominal control
level is achieved with a specified high probability. We refer to these
approaches as control-in-expectation (CiE) and control-in-probability (CiP),
respectively.

This article addresses outstanding issues in the EUM classification with two
objectives. The first is to resolve accuracy control for the prioritized
class. Since the nominal control level coincides with the empirical accuracy,
we conjecture that the observed under-coverage phenomenon is a manifestation
of the over-optimism bias well known in standard learning settings. We
develop asymptotic theory to support this conjecture and propose novel
threshold-adjustment techniques to refine the EUM classifier under both CiE
and CiP frameworks.
The second objective is to develop training-data-based methods for predicting
the performance of the resulting CiE and CiP classifiers. Reliable performance
estimation is essential for assessing potential deployment in future
applications. However, empirical estimation is generally unreliable, while
independent validation requires additional data. We address both
objectives within a unified framework by extending the cross-audit projection
method of \citet{huangc}, originally developed for risk prediction in
standard learning problems.

Throughout this article, the performance of an estimated classifier---namely
class-specific accuracies such as sensitivity and specificity---refers to
{\em prediction} performance, defined as the expected performance on independent
future data conditional on the training data. Performance evaluated on the
training data is explicitly referred to as {\em empirical} performance.

The remainder of the article is organized as follows. Section~2 develops
higher-order asymptotic theory to explain the under-coverage phenomenon of
the EUM classifier. Refinements for CiE and CiP frameworks are proposed in
Sections~3 and~4, respectively. Numerical studies are presented in Section~5,
and concluding remarks are provided in Section~6. Technical proofs are
deferred to the Appendix.

\section{Performance of the EUM classifier: an asymptotic analysis}

Let ${\bf M}_d$ denote a vector of $k \geq 2$ features under consideration
for class $d = 0, 1$. A linear classifier is specified by
$({\bf b}^\T, t)^\T$, where $\bf b$ and $t$ are combination coefficient vector
and threshold, respectively. We adopt the convention of assigning an
observation to class~0 if ${\bf b}^\T {\bf M}_d\leq t$, and to class~1
otherwise. The class-specific accuracies are given by
\[
\psi_0(t,{\bf b})=\Pr( {\bf b}^\T {\bf M}_0\leq t),\qquad
\psi_1(t,{\bf b})=\Pr( {\bf b}^\T {\bf M}_1> t).
\]
Since classification performance is scale-invariant in $({\bf b}^\T, t)^\T$,
we impose a normalization restriction $\|{\bf b}\|_1=1$, where
$\|\cdot\|_1$ denotes the $\ell_1$ norm. Class~0 is taken as the prioritized
class, and an oracle Neyman--Pearson classifier maximizes
$\psi_1(t,{\bf b})$ subject to maintaining $\psi_0(t,{\bf b})$ at a
prespecified level $\rho\in(0,1)$:
\begin{equation}
\max_{{\bf b}:\|{\bf b}\|_1=1,t}\psi_1(t,{\bf b})\qquad
\mbox{subject to }\psi_0(t,{\bf b})\geq \rho.
\label{oracle}
\end{equation}
As $\psi_0(t,{\bf b})$ is the cumulative distribution function of
${\bf b}^\T {\bf M}_0$, define its associated quantile function as
$\psi_0^{-1}(p,{\bf b})=\inf\{t:\psi_0(t,{\bf b})\geq p\}$. For a fixed
$\bf b$, the Neyman--Pearson classifier has threshold and class~1 accuracy
\begin{equation}
  \tau({\bf b})= \psi_0^{-1}(\rho,{\bf b}),\qquad
  \phi({\bf b})=\psi_1\{ \tau({\bf b}),{\bf b}\},
\end{equation}
respectively. The oracle classifier therefore uses the combination coefficient
\[
\bt=\arg\max_{{\bf b}:\|{\bf b}\|_1=1}\phi({\bf b})
\]
and threshold $\tau(\bt)$.

Consider a case-control study with $n_d$ independent replicates of ${\bf M}_d$,
${\bf M}_{d,[i]}$ for $i=1,\ldots,n_d$ and $d=0,1$, yielding a combined sample
of size $n=n_1+n_0$. To estimate the oracle classifier, \citet{huangb}
studied the EUM method:
\begin{equation}
\max_{{\bf b}:\|{\bf b}\|_1=1,t}\h{\psi}_1( t,{\bf b})\qquad
\mbox{subject to }\h{\psi}_0( t,{\bf b})\geq \rho,\label{geum}
\end{equation}
where $\h{\psi}_d(t,{\bf b})$ denotes the empirical counterpart of
$\psi_d(t,{\bf b})$ for $d=0,1$. The empirical analogues of
$\tau({\bf b})$ and $\phi({\bf b})$ are defined as
\begin{equation}
\h{\tau}({\bf b})=\h{\psi}_0^{-1}(\rho,{\bf b}),
\qquad \h{\phi}({\bf b})=\h{\psi}_1\{\h{\tau}({\bf b}),{\bf b}\},
\end{equation}
respectively. Profiling out $t$, estimation of the oracle combination
coefficient $\bt$ reduces to the optimization problem
\begin{equation}
  \max_{{\bf b}:\|{\bf b}\|_1=1}\h{\phi}({\bf b}),
  \label{eum}
\end{equation}
and a near maximizer $\hb$ satisfies
\begin{equation}
  \h{\phi}(\hb)\geq\max_{\bf b}\h{\phi}({\bf b})
  -o_p(n^{-2/3}),
  \qquad\|\hb\|_1=1.\label{hbdef}
\end{equation}
The threshold $\tau(\bt)$ is then estimated by $\h{\tau}(\hb)$.

\subsection{Effect of estimated combination in general}
Let $f_d(t,{\bf b})$ denote the probability density function of
${\bf b}^\T{\bf M}_d$, when it exists, for $d=0,1$. We impose the following
conditions.
\begin{condition}[Sample sizes]\label{con0}
As $n\rightarrow\infty$, $n_1/n_0$ converges to a finite positive constant.
\end{condition}
\begin{condition}[Identifiability and separation]\label{con1}
  Function $\phi({\bf b})$ and its maximizer $\bt$ satisfy
  $\phi(\bt)>\sup_{\|{\bf b}\|_1=1,\|{\bf b}-\bt\|_1\geq\epi}\phi({\bf b})$ for
  every $\epi>0$.
\end{condition}
\begin{condition}[Quantile smoothness of class~0]\label{con2}
  The density $f_0\{\tau(\bt),\bt\}$ exists and is strictly positive.
\end{condition}
\begin{condition}[Distribution smoothness]\label{con3} For $d=0,1$,
  there exists a $k\times(k-1)$ matrix ${\bf B}_d$ whose column vectors,
  together with $\bt$, are linearly independent, such that
  (i)~${\bf B}_d^\T{\bf M}_d$ is integrable; and (ii)~the conditional
  distribution function
  $\Pr({\bf b}^\T{\bf M}_d\leq t\mid {\bf B}_d^\T{\bf M}_d)$ has a bounded
  second derivative with respect to $t$, for $({\bf b}^\T,t)^\T$ in a
  neighborhood of $\{\bt^\T,\tau(\bt)\}^\T$.
\end{condition}

Consider a classifier with a general and consistent combination estimator
${\bf b}_n$, not necessarily $\hb$, together with threshold
$\h{\tau}({\bf b}_n)$. We characterize the bias in the class-specific
empirical accuracies.
\begin{prop}\label{decomp}
  Suppose that Conditions~\ref{con0}--\ref{con3} hold. For $d=0,1$ and any
  ${\bf b}_n=\bt+o_p(1)$,
    \begin{eqnarray}
      \lefteqn{\h{\psi}_d\{\h{\tau}({\bf b}_n),{\bf b}_n\}-\psi_d\{\h{\tau}({\bf b}_n),{\bf b}_n\}}\label{decom}\\
      &=& [\h{\psi}_d\{\tau(\bt),\bt\}-\psi_d\{\tau(\bt),\bt\}]\nonumber\\
      && \mbox{}+[\h{\psi}_d\{\tau({\bf b}_n),{\bf b}_n\}-\h{\psi}_d\{\tau(\bt),\bt\}
        -\psi_d\{\tau({\bf b}_n),{\bf b}_n\}+\psi_d\{\tau(\bt),\bt\}]\nonumber\\
      &&\mbox{}+O_p(n^{-3/4}\log^{3/4}n).\nonumber
    \end{eqnarray}
\end{prop}
In decomposition~(\ref{decom}), the first term has mean zero and is of order
$O_p(n^{-1/2})$. The second term, arising from the combination estimation, is
typically of an order between $O_p(n^{-1/2})$ and $O_p(n^{-3/4}\log^{3/4}n)$.
This term may therefore account for the bias in the empirical accuracies.
Meanwhile, the bias of class~0 control is directly linked to that in class~0
empirical accuracy since
$\h{\psi}_0\{\h{\tau}({\bf b}_n),{\bf b}_n\}=\rho+O(n^{-1})$. Notably,
the second term is not affected by the threshold estimation; the first term,
instead, is attributable to the estimated threshold.

\subsection{The case of the EUM classifier}
We now characterize the second term in decomposition~(\ref{decom}) for the
EUM classifier, that is, when ${\bf b}_n=\hb$. Note that this second term
reduces to $\h{\psi}_0\{\tau({\bf b}_n),{\bf b}_n\}-\h{\psi}_0\{\tau(\bt),\bt\}$
when $d=0$. Since classification is invariant under linear feature
transformation, we may assume without loss of generality that
$\bt=(1,{\bf 0}^\T)^\T$. To facilitate
the analysis, \citet{huangc} derived an approximation to
$\h{\phi}({\bf b})$ given by
\[
\xbar{\phi}({\bf b})=\lambda[\h{\psi}_0\{\tau({\bf b}),{\bf b}\}-\rho]+
\h{\psi}_1\{ \tau({\bf b}), {\bf b}\},
\]
where $\lambda=f_1\{\tau(\bt),\bt\}/f_0\{\tau(\bt),\bt\}$.
\begin{prop}\label{cra}
  Suppose that Conditions~\ref{con0}--\ref{con3} hold. Then, $\hb$ is
  consistent for $\bt$, and satisfies
  \[ \xbar{\phi}(\hb)\geq\max_{{\bf b}:{\bf b}_1=1}\xbar{\phi}({\bf b})-o_p(n^{-2/3}).
  \]
  Furthermore, as $n\rightarrow\infty$,
  \begin{equation}
    n^{2/3}\left[\begin{array}{c}\displaystyle
        \h{\psi}_0\{\tau({\bf g}_n), {\bf g}_n\}-\h{\psi}_0\{\tau(\bt),\bt\}\\
        \h{\psi}_1\{\tau({\bf g}_n),{\bf g}_n\}-\psi_1\{\tau({\bf g}_n), {\bf g}_n\}
        -\h{\psi}_1\{\tau(\bt),\bt\}+\psi_1\{\tau(\bt),\bt\}\\
        \xbar{\phi}({\bf g}_n)-\xbar{\phi}(\bt)
      \end{array}\right]\rightsquigarrow
    \left\{\begin{array}{c}\displaystyle
    W_0({\bf g})\\
    W_1({\bf g})\\
    Z({\bf g})\end{array}\right\},
  \end{equation}
  where ${\bf g}_n=(1,n^{-1/3}{\bf g})^\T$ for vector $\bf g$,
  $\rightsquigarrow$ denotes weak convergence, $W_0({\bf g})$ and
  $W_1({\bf g})$ are two independent mean-zero Gaussian processes, and
  $Z({\bf g})=\lambda W_0({\bf g})+W_1({\bf g})+{\bf g}^\T{\bf H}{\bf g}/2$
  for ${\bf H}=\nabla^2_{{\bf b}_{-1}}\phi({\bf b})\mid_{{\bf b}=\bt}$.
  Finally, $n^{2/3}[\h{\psi}_0\{\tau(\hb),\hb\}-\h{\psi}_0\{\tau(\bt),\bt\}]
  \rightsquigarrow W_0({\bf U})$ and
  $n^{2/3}[\h{\psi}_1\{\tau(\hb),\hb\}-\psi_1\{\tau(\hb), \hb\}
    -\h{\psi}_1\{\tau(\bt),\bt\}+\psi_1\{\tau(\bt),\bt\}]\rightsquigarrow W_1({\bf U})$, where
  ${\bf U}=\arg\max_{\bf g}Z({\bf g})$ and
  $E\{W_d({\bf U})\}>0$ for $d=0,1$.
\end{prop}

In light of the approximation $\xbar{\phi}({\bf b})$ to $\h{\phi}({\bf b})$,
it is not surprising that the second term in
decomposition~(\ref{decom})---which is of order
$O_p(n^{-2/3})$---exhibits an asymptotic positive bias, reflecting the
over-optimism of the EUM procedure. Accordingly, class~0 control
is second-order negatively biased relative to $\rho$. These results
corroborate the empirical results in Figure~\ref{fig1}, supporting our
conjecture that the under-coverage of the EUM classifier arises from
over-optimism.

\section{Refining the EUM classifier for control-in-expectation}

Given the estimated combination coefficient $\hb$, the corresponding
theoretical threshold for control level $\rho$ is $\tau(\hb)$. Consequently,
the empirical threshold $\h{\tau}(\hb)$ used in the EUM classifier is
expected to be negatively biased relative to $\tau(\hb)$.

\begin{theorem}\label{thm2}
  If Conditions~\ref{con0}--\ref{con3} hold, then
  \begin{eqnarray}
    \h{\tau}(\hb)-\tau(\hb)
    &=& f_0\{\tau(\bt),\bt\}^{-1}\left[
      \rho-\h{\psi}_0\{\tau(\bt),\bt\}\right.
      \label{bth}\\
      && \mbox{\hspace*{.2in}}+\left.
      \h{\psi}_0\{\tau(\bt),\bt\}
      -\h{\psi}_0\{\tau(\hb),\hb\}\right]
    +o_p(n^{-2/3}).\nonumber
  \end{eqnarray}
  Meanwhile, for any threshold $\h{t}=\tau(\hb)+O_p(n^{-1/2}\log^{1/2} n)$,
  \begin{eqnarray}    
    \psi_0(\h{t},\hb)-\rho
    &=&  f_0\{\tau(\bt),\bt\}\{\h{t}-\tau(\hb)\}
    +o_p(n^{-2/3}),\label{bsens}\\
    \psi_1(\h{t},\hb)-\phi(\hb) &=&
    - f_1\{\tau(\bt),\bt\}\{\h{t}-\tau(\hb)\}
    +o_p(n^{-2/3}).\label{bspec}
  \end{eqnarray}
\end{theorem}

The threshold $\h{\tau}(\hb)$ of the EUM classifier is consistent and
asymptotically first-order unbiased for $\tau(\hb)$.
However, it exhibits a negative second-order bias, since
$n^{2/3}[\h{\psi}_0\{\tau(\bt),\bt\}-\h{\psi}_0\{\tau(\hb),\hb\}]\rightsquigarrow-W_0({\bf U})$ as stated in Proposition~\ref{cra}. Equations~(\ref{bsens})
and~(\ref{bspec}) further characterize how bias of a consistent threshold
propagates to the class-specific accuracies. In particular, replacing
$\h{\tau}(\hb)$ with a second-order unbiased threshold yields a class~0
accuracy that is second-order unbiased for the control level $\rho$, while
the corresponding class~1 accuracy is second-order unbiased for $\phi(\hb)$.
These observations motivate a threshold-correction procedure.

\subsection{Threshold correction}
In light of equation~(\ref{bth}), the threshold behaves analogously to a
performance measure, with $\tau(\hb)$ serving as the prediction analogue
of the empirical quantity $\h{\tau}(\hb)$. This observation suggests adapting
risk prediction methods for threshold correction. Although $K$-fold
cross-validation (CV) is widely used for risk prediction, it may fail in
certain settings \citep{huangc}, particularly under cube-root asymptotics
\citep{kim}. We instead adapt the repeated $K$-fold cross-audit projection (CAP)
method of \citet{huangc}, focusing on the case $K=2$.

Let $r\geq1$ denote the number of CAP repetitions. We draw $r$ random
stratified half-and-half splits of the data by class, as evenly as possible
when $n_0$ or $n_1$ is odd. For the resulting $2r$ subsamples, we
compute the corresponding versions of $\h{\tau}(\hb)$ and average them to
obtain $\h{\tau}_{0.5}$. We also denote by $\h{\tau}_{0.5,\mathrm{cv}}$ the
corresponding repeated two-fold cross-validated threshold. Their difference
serves as an estimate of the bias of the empirical threshold at half the
sample size. Rescaling this bias estimate according to the $n^{-2/3}$ rate
yields the corrected threshold
\begin{equation}
  \h{\tau}_{\mathrm{c}}=\h{\tau}(\hb)
  -2^{-2/3}(\h{\tau}_{0.5}-\h{\tau}_{0.5,\mathrm{cv}}).
\end{equation}

\begin{prop}\label{prop1}
  Suppose that Conditions~\ref{con0}--\ref{con3} hold and
  the number $r$ of CAP repetitions is fixed. Then,
  $n^{2/3}\{\h{\tau}_{\mathrm{c}}-\tau(\hb)-f_0\{\tau(\bt),\bt\}^{-1}
  [\rho-\h{\psi}_0\{\tau(\bt),\bt\}]\}$ converges weakly to a mean-zero,
  nondegenerate distribution.
\end{prop}

The corrected threshold $\h{\tau}_{\mathrm{c}}$ is first-order asymptotically
equivalent to $\h{\tau}(\hb)$, but second-order asymptotically unbiased for
$\tau(\hb)$. The classifier $(\hb^\top,\h{\tau}_{\mathrm{c}})^\top$ is
referred to as cEUM. By Theorem~\ref{thm2}, the class~0 and class~1 accuracies
of the cEUM classifier are second-order asymptotically unbiased for $\rho$
and $\phi(\hb)$, respectively.

\subsection{Accuracy prediction and inference of the cEUM classifier}

The class-specific accuracies are given by
$\psi_{d,\mathrm{c}}\equiv\psi_d(\h{\tau}_{\mathrm{c}},\hb)$ for $d=0,1$. We
aim to estimate these quantities using only the training data. A natural
estimator for $\psi_{0,\mathrm{c}}$ is $\rho$, which is second-order unbiased.
Estimating $\psi_{1,\mathrm{c}}$ is more challenging because the threshold
$\h{\tau}_{\mathrm{c}}$ is obtained via a CAP procedure. A straightforward
nested CAP implementation would be computationally prohibitive. We therefore
proceed indirectly by targeting $\phi(\hb)$ instead. As shown in
equation~(\ref{bspec}), any estimator that is second-order unbiased for
$\phi(\hb)$ is also second-order unbiased for
$\psi_1(\h{\tau}_{\mathrm{c}},\hb)$. A standard CAP procedure can therefore
be used for this purpose.

We adopt a repeated two-fold CAP procedure with $r$ repetitions, using the
same resampling scheme as above. For the resulting $2r$ subsamples, we
compute the corresponding versions of $\h{\phi}(\hb)$ and average them to
obtain $\h{\phi}_{0.5}$. The repeated two-fold cross-validated estimator
$\h{\phi}_{0.5,\mathrm{cv}}$ is obtained by learning the combination on one
subsample and then estimating the threshold, followed by evaluating
class~1 accuracy on the other subsample within each split.
The CAP estimator is defined as
\begin{equation}
\h{\phi}_\mathrm{cap}=Q^{-1}\left[Q\{\h{\phi}(\hb)\}
-2^{-2/3}\{Q(\h{\phi}_{0.5})-Q(\h{\phi}_{0.5,\mathrm{cv}})\}\right],\label{cpp}
\end{equation}
where $Q(\cdot)$ is a monotone transformation. For range preservation, we took
$Q(\cdot)$ to be the standard normal quantile function in our numerical studies.

\begin{coroll}\label{coro1}
  Suppose that Conditions~\ref{con0}--\ref{con3} hold,
  the number $r$ of CAP repetitions is fixed, and the monotone transformation
  $Q(\cdot)$
  is differentiable at $\phi(\bt)$ with a nonzero derivative. Then,
  $n^{2/3}[\h{\phi}_\mathrm{cap}-\psi_1(\h{\tau}_{\mathrm{c}},\hb)-\h{\psi}_1\{\tau(\bt),\bt\}+\psi_1\{\tau(\bt),\bt\}]$
  converges weakly to a mean-zero, nondegenerate distribution.
\end{coroll}

Let $\mbox{AN}$ denote asymptotically normal. First-order asymptotics theory
yields
\[
\left(\begin{array}{c}\h{\psi}_{0,\mathrm{c}}\\ \h{\psi}_{1,\mathrm{c}}
\end{array}\right)\equiv
\left(\begin{array}{c}\rho\\ \h{\phi}_\mathrm{cap}
\end{array}\right)\sim
\mbox{AN}\left\{\left(\begin{array}{c}\psi_{0,\mathrm{c}}\\ \psi_{1,\mathrm{c}}
\end{array}\right),
\left(\begin{array}{cc}\psi_{0,\mathrm{c}}(1-\psi_{0,\mathrm{c}})/n_0 & 0\\
    0 & \psi_{1,\mathrm{c}}(1-\psi_{1,\mathrm{c}})/n_1\end{array}\right)\right\},
\]
which provides the basis for inference on $\psi_{0,\mathrm{c}}$
and $\psi_{1,\mathrm{c}}$. The procedure follows that of
\citet[section~3.3]{huangc}. Specifically, for each class $d=0,1$, a
lower accuracy bound with nominal confidence level $\alpha_d$ is obtained as
the smaller solution $q$ to
\[
\frac{n_d(q-\h{\psi}_{d,\mathrm{c}})^2}{q(1-q)}=z_{\alpha_d}^2,
\]
where $z_{\alpha_d}$ denotes the $\alpha_d$ quantile of the standard normal
distribution. Because of the asymptotic independence, the joint confidence
level of the two class-specific accuracy bounds is asymptotically
$\alpha_0\alpha_1$.

\section{Classification with control-in-probability}

Control-in-probability (CiP) provides an alternative strategy that aims to
maintain the class~0 accuracy at the control level $\rho$ with a prespecified,
typically high, probability $\delta$. We extend the EUM framework to
accommodate this criterion. Let $\rho_n$ denote a sequence of control levels converging to $\rho$, and define
$\h{\tau}({\bf b};\rho_n)=\h{\psi}_0^{-1}(\rho_n,{\bf b})$.
Equation~(\ref{bth}) extends to this more general setting.
\begin{coroll}\label{coro41}
  If Conditions~\ref{con0}--\ref{con3} hold and
  $\rho_n=\rho+O\{n^{-1/2}(\log n)^{1/2}\}$,
  \begin{eqnarray}
    \h{\tau}(\hb;\rho_n)-\tau(\hb)
    &=& f_0\{\tau(\bt),\bt\}^{-1}\left[
      \rho_n-\h{F}_0\{\tau(\bt),\bt\}\right.
      \label{bth2}\\
      && \mbox{\hspace*{.2in}}+\left.
      \h{F}_0\{\tau(\bt),\bt\}
      -\h{F}_0\{\tau(\hb),\hb\}\right]
    +o_p(n^{-2/3}).\nonumber
  \end{eqnarray}
\end{coroll}

Following equations~(\ref{bth2}) and~(\ref{bsens}),
\begin{equation}
  \psi_0\{\h{\tau}(\hb;\rho_n),\hb\}-\rho
  =\rho_n-\h{F}_0\{\tau(\bt),\bt\}+O_p(n^{-2/3}).\label{bsens2}
\end{equation}
This result suggests that asymptotic control-in-probability can be achieved
by setting $\rho_n=q_\mathrm{Binom}(\delta;n_0,\rho)/n_0$,
where $q_\mathrm{Binom}(\delta;n_0,\rho)$ denotes the $\delta$-quantile of the
binomial distribution with size $n_0$ and success probability $\rho$. The
resulting classifier is referred to as EUM$_\mathrm{p}$.
Let $\h{\tau}_{\mathrm{c}}(\rho_n)$ denote the corrected threshold obtained by
applying the procedure of Section~3.1 with $\rho_n$ replacing $\rho$. The
corresponding classifier, which uses threshold $\h{\tau}_{\mathrm{c}}(\rho_n)$,
is referred to as cEUM$_\mathrm{p}$.

\begin{remark}\rm
  Both EUM$_\mathrm{p}$ and cEUM$_\mathrm{p}$ take the combination component
  $\hb$ of the EUM classifier. It can be shown that $\hb$ is a near maximizer of
  $\h{\phi}({\bf b};\rho_n)\equiv\h{\psi}_1\{\h{\tau}({\bf b};\rho_n),{\bf b}\}$
  subject to $\|{\bf b}\|_1=1$.
  Consequently, EUM$_\mathrm{p}$ and cEUM$_\mathrm{p}$ can be equivalently
  obtained by applying the EUM and cEUM
  procedures, respectively, with $\rho$ replaced by $\rho_n$.
\end{remark}

\begin{coroll}\label{thmcip}
  Suppose that Conditions~\ref{con0}--\ref{con3} hold and, if applicable,
  the number $r$ of CAP repetitions is fixed. For
  $\h{t}=\h{\tau}(\hb;\rho_n)$ or $\h{\tau}_\mathrm{c}(\rho_n)$ where
  $\rho_n=q_\mathrm{Binom}(\delta;n_0,\rho)/n_0$, the
  thresholds of EUM$_\mathrm{p}$ or cEUM$_\mathrm{p}$, respectively,
  $\Pr\{\psi_0(\h{t};\hb)\geq\rho\}\rightarrow\delta$.
\end{coroll}

For both the EUM$_\mathrm{p}$ and cEUM$_\mathrm{p}$ classifiers, the probability
that the class~0 accuracy exceeds the control level converges to the
prespecified value $\delta$. Our simulations indicate a clear advantage of
cEUM$_\mathrm{p}$ over EUM$_\mathrm{p}$ in finite samples, although a formal
theoretical guarantee remains to be established. Nevertheless, in terms of
class~0 accuracy relative to $\rho_n$, it can be shown that the
cEUM$_\mathrm{p}$ classifier is second-order unbiased, but the
EUM$_\mathrm{p}$ classifier is not. 

The EUM$_\mathrm{p}$ and cEUM$_\mathrm{p}$ classifiers are approximately
equivalent to the EUM and cEUM classifiers, respectively, when $\delta=0.5$.
This equivalence does not hold in general, even asymptotically. In practice,
$\delta>0.5$ is typically chosen, in which case EUM$_\mathrm{p}$ and
cEUM$_\mathrm{p}$ impose a more stringent class~0 accuracy control than their
CiE counterparts under the same nominal control level $\rho$.

For the cEUM$_\mathrm{p}$ classifier, we further consider performance
estimation. Since class~0 accuracy is now controlled in probability, class~1
accuracy becomes the only quantity of interest. The estimation and inference
procedures in Section~3.2 can be directly applied with $\rho_n$ replacing
$\rho$, yielding an asymptotically second-order unbiased point estimator as
well as an associated accuracy bound.

\section{Numerical studies}

Simulations were conducted to evaluate the proposed Neyman--Pearson classifiers
and their associated performance prediction procedures. These methods were
also illustrated using a breast cancer detection study. In all numerical
experiments, the performance metric of interest was specificity at a
controlled sensitivity level. For the control-in-expectation framework, this
corresponded to 95\% sensitivity ($\rho=0.95$), while for the
control-in-probability framework it corresponded to 90\% sensitivity with
probability 0.9 ($\rho=0.9$, $\delta=0.9$).
The algorithm of \cite{huangb} was used to compute the estimated linear
combination coefficients. For the adapted CAP procedures, the number of
repetitions was set to $r=16$.

\subsection{Simulations}

\begin{table}[t]
  \def~{\hphantom{0}}
  \caption{\label{tab2}Simulation results on the performance of
    the EUM and cEUM
    classifiers for maximum specificity at controlled 95\% sensitivity}
\centerline{
\begin{tabular}{
    c@{\hspace*{4pt}}c@{\hspace*{8pt}}
    c@{\hspace*{4pt}}c@{\hspace*{-1pt}}c
    c@{\hspace*{4pt}}c@{\hspace*{-1pt}}c
    c@{\hspace*{4pt}}c@{\hspace*{4pt}}c@{\hspace*{8pt}}c
    c@{\hspace*{4pt}}c@{\hspace*{-1pt}}c
    c@{\hspace*{4pt}}c@{\hspace*{-1pt}}c
    c@{\hspace*{4pt}}c@{\hspace*{4pt}}c
  }\\[-10pt]\hline\\[-12pt]
$n_0,$ & $k$ & \multicolumn{2}{l}{EUM} && \multicolumn{6}{l}{cEUM}
            && \multicolumn{2}{l}{EUM} && \multicolumn{6}{l}{cEUM}\\
  \cline{3-4}\cline{6-11}\cline{13-14}\cline{16-21}
  $\phantom{,}n_1$&& \multicolumn{1}{l}{SN} & \multicolumn{1}{l}{SP} && \multicolumn{2}{l}{SN} && \multicolumn{3}{l}{SP} &&
  \multicolumn{1}{l}{SN} & \multicolumn{1}{l}{SP} && \multicolumn{2}{l}{SN} && \multicolumn{3}{l}{SP}\\
  && \multicolumn{1}{l}{M} & \multicolumn{1}{l}{M} && \multicolumn{1}{l}{M} & \multicolumn{1}{l}{C} && \multicolumn{1}{l}{M} & \multicolumn{1}{l}{E} & \multicolumn{1}{l}{C} &&
  \multicolumn{1}{l}{M} & \multicolumn{1}{l}{M} && \multicolumn{1}{l}{M} & \multicolumn{1}{l}{C} && \multicolumn{1}{l}{M} & \multicolumn{1}{l}{E} & \multicolumn{1}{l}{C}\\
\cline{3-11}\cline{13-21}
&& \multicolumn{9}{c}{Scenario A} && \multicolumn{9}{c}{Scenario B}\\
100 & 3 &91.9 & 52.7 && 94.2 & 94.4 && 45.3 & 45.4 & 94.7
       &&91.8 & 50.2 && 94.2 & 93.8 && 43.9 & 43.7 & 96.2\\
    & 6 &89.4 & 55.9 && 94.0 & 94.3 && 43.5 & 43.2 & 94.4
       &&89.4 & 52.8 && 94.0 & 94.3 && 41.9 & 41.2 & 95.2\\
200 & 3 &93.0 & 50.8 && 94.4 & 94.0 && 46.1 & 46.5 & 94.1
       &&93.0 & 48.6 && 94.4 & 93.9 && 44.6 & 44.8 & 94.7\\
    & 6 &91.6 & 53.3 && 94.2 & 91.1 && 45.4 & 45.6 & 95.4
       &&91.5 & 50.6 && 94.2 & 90.2 && 43.8 & 43.8 & 95.2\\
300 & 3 &93.5 & 50.0 && 94.6 & 93.5 && 46.1 & 45.8 & 96.9
       &&93.5 & 47.8 && 94.6 & 93.7 && 44.4 & 44.2 & 95.6\\
    & 6 &92.5 & 51.6 && 94.7 & 95.6 && 44.6 & 44.3 & 95.6
       &&92.6 & 49.1 && 94.8 & 96.3 && 43.2 & 42.7 & 96.9\\
400 & 3 &93.9 & 49.1 && 94.8 & 94.0 && 45.9 & 45.9 & 95.6
       &&93.9 & 47.1 && 94.8 & 95.5 && 44.4 & 44.3 & 94.9\\
    & 6 &92.9 & 51.0 && 94.6 & 93.5 && 45.4 & 45.4 & 94.6
       &&92.9 & 48.6 && 94.6 & 93.1 && 44.0 & 44.0 & 94.7\\
500 & 3 &94.0 & 48.9 && 94.8 & 94.1 && 45.9 & 45.7 & 95.6
       &&93.9 & 47.1 && 94.8 & 93.4 && 44.5 & 44.4 & 94.8\\
    & 6 &93.2 & 50.5 && 94.8 & 93.9 && 45.2 & 45.1 & 94.4
       &&93.2 & 48.3 && 94.7 & 93.7 && 44.0 & 43.8 & 95.0\\
\cline{3-11}\cline{13-21}
&& \multicolumn{9}{c}{Scenario C} && \multicolumn{9}{c}{Scenario D}\\
100 & 3 &92.2 & 48.1 && 94.1 & 94.3 && 42.8 & 42.4 & 95.4
       &&91.4 & 51.7 && 94.1 & 94.5 && 43.0 & 42.9 & 95.9\\
    & 6 &90.0 & 50.6 && 94.1 & 95.1 && 40.4 & 39.1 & 96.4
       &&89.0 & 55.0 && 93.9 & 93.9 && 41.4 & 40.9 & 96.5\\
200 & 3 &93.3 & 46.7 && 94.4 & 93.6 && 43.4 & 43.3 & 94.0
       &&92.9 & 48.8 && 94.4 & 91.6 && 43.4 & 43.7 & 94.3\\
    & 6 &91.9 & 48.5 && 94.3 & 91.6 && 42.1 & 41.6 & 94.5
       &&91.4 & 51.6 && 94.2 & 91.4 && 42.9 & 42.6 & 95.5\\
300 & 3 &93.7 & 46.1 && 94.7 & 93.3 && 43.5 & 43.1 & 96.6
       &&93.5 & 47.7 && 94.8 & 95.3 && 43.1 & 42.9 & 94.9\\
    & 6 &92.7 & 47.4 && 94.6 & 94.3 && 42.2 & 41.5 & 96.5
       &&92.3 & 50.0 && 94.6 & 94.0 && 42.5 & 42.3 & 95.2\\
400 & 3 &94.0 & 45.7 && 94.7 & 92.9 && 43.5 & 43.6 & 93.5
       &&93.7 & 47.1 && 94.7 & 95.2 && 43.4 & 43.5 & 93.4\\
    & 6 &93.1 & 47.0 && 94.6 & 92.8 && 42.8 & 42.4 & 96.0
       &&92.8 & 48.9 && 94.6 & 94.1 && 42.8 & 42.9 & 94.0\\
500 & 3 &94.2 & 45.4 && 94.8 & 94.3 && 43.5 & 43.3 & 95.5
       &&93.9 & 46.7 && 94.9 & 95.3 && 43.2 & 43.2 & 95.7\\
    & 6 &93.4 & 46.4 && 94.8 & 93.5 && 42.5 & 42.0 & 96.6
       &&93.1 & 48.3 && 94.8 & 93.4 && 42.6 & 42.5 & 95.8\\
\hline\\[-8pt]
\multicolumn{21}{l}{\footnotesize
  SN: sensitivity, SP: specificity.}\\
\multicolumn{21}{l}{\footnotesize
  M: mean ($\times100$), E: mean of predicted performance
  ($\times100$), C: coverage of 95\% accuracy bound.}
\end{tabular}}
\end{table}

The same set-ups as in \cite{huangb} were adopted, mimicking biomarker-based
cancer detection. All control biomarkers were
independent and identically distributed as standard normal, whereas the case
biomarkers varied across simulation scenarios. For
set-ups with $k=3$ biomarkers, four case distributions were considered:
\begin{list}{}{\itemsep 0ex\parsep 0ex\topsep 0ex}
\item[{\it Scenario~A.}] Independent and identically distributed normal
  variables with mean 0.9 and variance 1;
\item[{\it Scenario~B.}] Independent normal variables with common
  mean of 0.8 but heterogeneous variances of 0.5, 1, and 2;
\item[{\it Scenario~C.}] Jointly normal variables with common mean 1,
  variances of 0.5, 1, and 2, and pairwise correlation 0.5;
\item[{\it Scenario~D.}] A mixture of two distributions with three
  independent normal biomarkers: with probability
  2/3, means (1.7, 1.7, 0) and variances (0.5, 2, 1); and with probability 1/3,
  means (0, 0, 1.7) and common variance of 1.
\end{list}
In the first three scenarios, all case biomarkers were elevated relative to
their control counterparts, but the scenarios differed in dependence structure
and variability. Scenario~D captured cancer heterogeneity with two subtypes,
responsible for the elevation of the first two biomarkers and the last one
separately. For each scenario, we additionally considered a setting with
$k=6$ biomarkers by augmenting the model with three independent,
non-informative standard normal case biomarkers. At the 95\% sensitivity control
level, the oracle specificities were 0.466, 0.452, 0.444, and 0.442 under
Scenarios~A, B, C, and D, respectively. Case and control sizes were taken
to be equal, $n_0=n_1$, ranging from 100 to 500. All results were based on
1,000 simulation replications for each setting.

\begin{table}[t!]
  \def~{\hphantom{0}}
  \caption{\label{tab3}Simulation results on the performance of estimated
    classifiers targeting controlled 90\% sensitivity with
    probability 0.9}
\centerline{
\begin{tabular}{
    c@{\hspace*{3pt}}c@{\hspace*{8pt}}
    c@{\hspace*{3pt}}c@{}c@{\hspace*{4pt}}
    c@{\hspace*{3pt}}c@{}c@{\hspace*{4pt}}
    c@{\hspace*{3pt}}c@{\hspace*{3pt}}c@{\hspace*{3pt}}c@{}
    c@{\hspace*{6pt}}
    c@{\hspace*{3pt}}c@{}c@{\hspace*{4pt}}
    c@{\hspace*{3pt}}c@{}c@{\hspace*{4pt}}
    c@{\hspace*{3pt}}c@{\hspace*{3pt}}c@{\hspace*{3pt}}c}\\[-10pt]\hline\\[-12pt]
  $n_0,$ & $k$ &  \multicolumn{2}{l}{eLDA}
               && \multicolumn{2}{l}{EUM$_\mathrm{p}$}
               && \multicolumn{4}{l}{cEUM$_\mathrm{p}$}
               && \multicolumn{2}{l}{eLDA}
               && \multicolumn{2}{l}{EUM$_\mathrm{p}$}
               && \multicolumn{4}{l}{cEUM$_\mathrm{p}$}\\
\cline{3-4}\cline{6-7}\cline{9-12}
\cline{14-15}\cline{17-18}\cline{20-23}
  $\phantom{,}n_1$&& \multicolumn{1}{l}{SN} & \multicolumn{1}{l}{SP}
  && \multicolumn{1}{l}{SN} & \multicolumn{1}{l}{SP}
  && \multicolumn{1}{l}{SN} & \multicolumn{3}{l}{SP}
  && \multicolumn{1}{l}{SN} & \multicolumn{1}{l}{SP}
  && \multicolumn{1}{l}{SN} & \multicolumn{1}{l}{SP}
  && \multicolumn{1}{l}{SN} & \multicolumn{3}{l}{SP}\\
  && \multicolumn{1}{l}{P} & \multicolumn{1}{l}{M}
   && \multicolumn{1}{l}{P} & \multicolumn{1}{l}{M}
   && \multicolumn{1}{l}{P} & \multicolumn{1}{l}{M}
    & \multicolumn{1}{l}{E} & \multicolumn{1}{l}{C}
   && \multicolumn{1}{l}{P} & \multicolumn{1}{l}{M}
   && \multicolumn{1}{l}{P} & \multicolumn{1}{l}{M}
   && \multicolumn{1}{l}{P} & \multicolumn{1}{l}{M}
    & \multicolumn{1}{l}{E} & \multicolumn{1}{l}{C}\\
  \cline{3-12}\cline{14-23}
  && \multicolumn{10}{c}{Scenario A} && \multicolumn{10}{c}{Scenario B}\\
100 & 3 &80.8 & 56.1 && 82.3 & 51.7 && 89.9 & 48.5 & 48.4 & 94.7
       &&87.4 & 49.3 && 80.0 & 48.9 && 88.6 & 45.8 & 45.6 & 95.2\\
    & 6 &76.7 & 55.4 && 69.0 & 52.1 && 88.8 & 46.0 & 45.1 & 96.3
       &&83.8 & 48.7 && 66.9 & 49.0 && 90.0 & 43.3 & 42.2 & 96.7\\
200 & 3 &80.4 & 57.9 && 73.8 & 55.8 && 84.7 & 53.7 & 53.7 & 94.1
       &&91.7 & 51.0 && 74.3 & 52.4 && 85.5 & 50.5 & 50.7 & 93.7\\
    & 6 &77.0 & 57.5 && 64.3 & 56.1 && 83.4 & 52.4 & 52.4 & 95.7
       &&88.6 & 50.6 && 60.7 & 53.1 && 83.7 & 49.6 & 49.4 & 94.9\\
300 & 3 &79.6 & 58.4 && 82.9 & 55.9 && 89.7 & 54.4 & 54.3 & 95.7
       &&92.2 & 51.6 && 81.0 & 52.6 && 88.8 & 51.2 & 51.0 & 95.9\\
    & 6 &82.0 & 57.9 && 75.0 & 55.8 && 92.5 & 53.0 & 52.8 & 95.1
       &&93.2 & 51.1 && 72.8 & 52.7 && 90.6 & 50.0 & 49.6 & 96.3\\
400 & 3 &80.4 & 58.7 && 86.1 & 56.4 && 91.8 & 55.0 & 54.8 & 95.6
       &&93.9 & 51.8 && 84.6 & 53.1 && 92.2 & 51.8 & 51.6 & 96.2\\
    & 6 &76.7 & 58.6 && 72.8 & 56.8 && 89.2 & 54.3 & 54.0 & 95.0
       &&92.3 & 51.7 && 73.0 & 53.4 && 89.3 & 51.0 & 50.7 & 96.0\\
500 & 3 &78.1 & 59.1 && 83.1 & 57.1 && 90.4 & 56.0 & 55.9 & 94.9
       &&94.4 & 52.2 && 82.7 & 53.7 && 90.2 & 52.7 & 52.6 & 95.3\\
    & 6 &77.7 & 58.9 && 74.9 & 57.2 && 89.9 & 55.1 & 54.9 & 94.2
       &&93.5 & 52.0 && 73.5 & 53.9 && 89.2 & 52.0 & 51.9 & 95.0\\
  \cline{3-12}\cline{14-23}
  && \multicolumn{10}{c}{Scenario C} && \multicolumn{10}{c}{Scenario D}\\
100 & 3 &51.8 & 51.7 && 82.0 & 46.8 && 88.3 & 44.3 & 44.0 & 94.7
       &&54.9 & 56.2 && 79.5 & 49.6 && 88.8 & 45.8 & 45.7 & 95.3\\
    & 6 &50.7 & 50.6 && 71.6 & 46.5 && 90.5 & 40.9 & 39.5 & 96.6
       &&56.0 & 55.2 && 64.8 & 50.1 && 88.3 & 43.3 & 42.5 & 96.7\\
200 & 3 &37.2 & 53.6 && 76.8 & 50.7 && 85.2 & 49.1 & 48.9 & 95.4
       &&49.9 & 57.5 && 74.0 & 53.5 && 83.8 & 51.0 & 51.1 & 96.1\\
    & 6 &38.3 & 53.0 && 64.2 & 50.9 && 83.1 & 47.4 & 47.0 & 94.7
       &&49.7 & 57.0 && 60.5 & 54.2 && 84.2 & 49.9 & 49.5 & 95.8\\
300 & 3 &27.3 & 54.6 && 84.0 & 50.8 && 88.5 & 49.7 & 49.5 & 95.6
       &&45.0 & 58.2 && 82.4 & 53.6 && 89.7 & 51.7 & 51.6 & 94.7\\
    & 6 &30.7 & 53.9 && 72.8 & 50.8 && 89.1 & 48.2 & 47.8 & 95.5
       &&45.8 & 57.8 && 70.9 & 54.1 && 88.8 & 50.9 & 50.7 & 95.7\\
400 & 3 &20.5 & 54.9 && 81.6 & 51.7 && 88.5 & 50.6 & 50.6 & 94.3
       &&39.2 & 58.7 && 83.6 & 54.4 && 91.4 & 52.7 & 52.7 & 95.0\\
    & 6 &20.7 & 54.6 && 73.6 & 51.7 && 88.8 & 49.4 & 49.0 & 96.2
       &&38.5 & 58.3 && 74.2 & 54.7 && 90.0 & 51.8 & 51.7 & 95.3\\
500 & 3 &16.6 & 55.1 && 84.3 & 52.0 && 89.8 & 51.1 & 51.0 & 93.7
       &&34.7 & 58.9 && 83.0 & 55.0 && 91.2 & 53.6 & 53.6 & 95.1\\
    & 6 &17.5 & 54.8 && 77.2 & 51.9 && 89.6 & 50.0 & 49.7 & 95.7
       &&34.1 & 58.6 && 74.1 & 55.2 && 89.4 & 52.8 & 52.8 & 94.1\\
\hline\\[-8pt]
\multicolumn{23}{l}{\footnotesize
  SN: sensitivity, SP: specificity.}\\
\multicolumn{23}{l}{\footnotesize
  P: probability (\%) of exceeding the
  control level, M: mean ($\times100$),}\\
\multicolumn{23}{l}{\footnotesize
   E: mean of estimated
  specificity ($\times100$), C: coverage of 95\% accuracy bound.}
\end{tabular}}
\end{table}

We begin with the EUM and cEUM classifiers for control-in-expectation framework.
Figure~\ref{fig1}, presented in Section~1, shows the class-specific
accuracies of the EUM classifier under Scenario~A; results for Scenarios~B--D
exhibited similar patterns and are omitted for brevity. Table~\ref{tab2}
summarizes the accuracies of the EUM and cEUM classifiers, along with the
proposed performance prediction results for the cEUM classifier.
The cEUM classifier achieved sensitivity close to the nominal control
level of 0.95 on average, representing a substantial improvement over the
EUM classifier. Moreover, the proposed prediction method accurately tracked
the performance of the cEUM classifier on average, with coverage
probabilities of the accuracy bounds close to the nominal level.

We next consider the estimated classifiers under the control-in-probability
framework. In addition to the proposed EUM$_\mathrm{p}$ and cEUM$_\mathrm{p}$
classifiers, we included the eLDA classifier of \cite{wang}, which is
developed under the linear discriminant analysis model, for comparison.
Table~\ref{tab3} summarizes the simulation results.
Among the four scenarios considered, the linear discriminant analysis model
is correctly specified only in Scenario~A. The eLDA classifier fell
substantially short of the nominal 90\% probability level for achieving the
desired sensitivity control even under Scenario A, and performed worse under
Scenarios~C and~D. The EUM$_\mathrm{p}$ classifier also failed to attain the
nominal probability level, although its performance was more stable across
the four scenarios. In contrast, the cEUM$_\mathrm{p}$ classifier performed
best overall, with the probability of achieving the desired sensitivity
level close to the nominal target. In addition, its estimated specificity
had little bias, and the corresponding accuracy
bound exhibited satisfactory coverage.

\subsection{Application to breast cancer detection}

To illustrate the proposed methods, we applied them to a breast cancer
detection study using demographic characteristics and blood measurements.
The Breast Cancer Coimbra dataset \citep{patricio} consists of 64 women with
breast cancer and 52 healthy controls. The features used in the analysis were
glucose, resistin, age, and BMI, all of which were log-transformed. The
cancer class was taken as the prioritized class for accuracy control.

\begin{table}[t]
  \caption{\label{tab4}Breast cancer detection: classification under
    Neyman--Pearson paradigm}
\centerline{
  \begin{tabular}{lrc lcc lcc}\\[-10pt]\hline\\[-12pt]
    \multicolumn{9}{l}{\parbox{4.4in}{Performance metric: specificity at 95\% sensitivity}}\\
    \multicolumn{9}{l}{control of 95\% sensitivity: in expectation}\\[3pt]
    & \multicolumn{1}{c}{coef} &&& \multicolumn{1}{c}{thresh}
    && \multicolumn{1}{c}{PE} & \multicolumn{1}{c}{AB}\\
log(glucose)  &    0.623 && EUM  & 2.290\\
log(resistin) &    0.154\\
log(age)      & $-0.139$ && cEUM & 2.259 & SN & 0.950 & 0.884\\
log(BMI)      & $-0.093$ &&      &       & SP & 0.376 & 0.274\\[1pt]\hline\\[-12pt]
    \multicolumn{9}{l}{Performance metric: specificity at 90\% sensitivity}\\
    \multicolumn{9}{l}{control of 90\% sensitivity: in probability of 0.9}\\[3pt]
    & \multicolumn{1}{c}{coef} &&& \multicolumn{1}{c}{thresh}
    && \multicolumn{1}{c}{PE} & \multicolumn{1}{c}{AB}\\
log(glucose)  &    0.677 && EUM$_\mathrm{p}$  & 2.441\\
log(resistin) &    0.106\\
log(age)      & $-0.117$ && cEUM$_\mathrm{p}$ & 2.433\\
log(BMI)      & $-0.100$ &&      &       & SP & 0.351 & 0.252\\[1pt]\hline\\[-8pt]
\multicolumn{9}{l}{\footnotesize
  coef: combination coefficient, thresh: threshold, PE: point estimate,}\\
\multicolumn{9}{l}{\footnotesize
  AB: 95\% accuracy bound.}\\
\multicolumn{9}{l}{\footnotesize
  SN: sensitivity, SP: specificity.}
\end{tabular}}
\end{table}

The results are reported in Table~\ref{tab4}. Under the control-in-expectation
framework with 95\% sensitivity, the EUM and cEUM classifiers were obtained.
These classifiers shared the same combination coefficients, but the cEUM
classifier employed a smaller threshold to correct the under-control
phenomenon of EUM. In addition, estimated accuracies and 95\% accuracy bounds
are reported for the cEUM classifier.
We also report the EUM$_\mathrm{p}$ and cEUM$_\mathrm{p}$ classifiers under
the control-in-probability framework targeting 90\% sensitivity with
probability 0.9, along with the estimated specificity and 95\% specificity
bounds for the cEUM$_\mathrm{p}$ classifier.

\section{Discussion}

Neyman--Pearson classification is complicated by the fact that performance
determines the acceptability of a classifier. As a consequence, the
under-coverage phenomenon of the EUM classifier is closely linked to its
over-optimism, a connection that not only provides insight into the behavior
of the EUM procedure but also motivates the refined classifiers developed
in this work. Performance prediction and inference for these methods are
developed in parallel within a unified framework.

The theoretical development of this article is based on higher-order
asymptotics under mild regularity conditions, in contrast to the nonasymptotic
learning theory of \cite{cannon} and \cite{scott05}. While the latter offers
greater generality, its resulting bounds are often highly conservative and
therefore of limited practical relevance.

Several directions warrant further investigation. First, the proposed
performance inference procedures are justified via first-order asymptotics
and may benefit from higher-order refinements. Second, the current inference
framework is limited to a single refined EUM classifier, and extending it to
comparisons with competing classifiers would be of practical interest.
Finally, incorporating feature selection---particularly in the presence of
non-informative features---may further improve performance and remains an
active area for future research.


\section*{Appendix: Technical details}

The asymptotic analyses extend those in \cite{huangb} to address estimated
classification. The regularity conditions, Conditions~\ref{con0}--\ref{con3},
are adapted from those in \cite{huangb}. Although Condition~\ref{con3}
is more general in allowing for linear feature transformations, it still
guarantees smoothness of the class-specific marginal distributions around
$\tau(\bt)$ for combinations with coefficients in a neighborhood of $\bt$.
Consequently, the results and arguments of \cite{huangb}
continue to apply. At the same time, some conditions are strengthened
for simplicity. In particular, Condition~\ref{con2} together with
part~(ii) of Condition~\ref{con3} in the case of $d=0$ imply a version of
\citet[condition~6(ii)]{huangb}.

\subsection*{Proof of Proposition~\ref{decomp}}
  Since $\h{\tau}(c{\bf b})=c\h{\tau}({\bf b})$ and
  $\tau(c{\bf b})=c\tau({\bf b})$ for all $c>0$ and $\bf b$,
  \citet[lemma~3.2]{huangb} extends to yield
  \[ \h{\tau}({\bf b}_n)=\tau({\bf b}_n)+O_p(n^{-1/2}\log^{1/2} n). \]
  Likewise, \citet[lemma~3.3]{huangb} extends to accommodate ${\bf b}_n$.
  Together, they imply the assertion.

\subsection*{Proof of Proposition~\ref{cra}}
  The first half of the proposition follows the arguments used in the proof of
  \citet[theorems~3.1 and~3.5]{huangb}. The weak convergence of
  $n^{2/3}[\h{\psi}_0\{\tau(\hb),\hb\}-\h{\psi}_0\{\tau(\bt),\bt\}]$ and
  $n^{2/3}[\h{\psi}_1\{\tau(\hb),\hb\}-\psi_1\{\tau(\hb), \hb\}
    -\h{\psi}_1\{\tau(\bt),\bt\}+\psi_1\{\tau(\bt),\bt\}]$, together with the
  existence and positivity of $E\{W_d({\bf U})\}$ for $d=0,1$,
  can be established by adapting the proof of \citet[theorem~3]{huangc}.

\subsection*{Proof of Theorem~\ref{thm2}}
  By Proposition~\ref{cra},
  \begin{eqnarray*}
    f_0\{\tau(\hb),\hb\} &=& f_0\{\tau(\bt),\bt\}+O_p(n^{-1/3}),\\
    \rho-\h{\psi}_0\{\tau(\hb),\hb\} &=& [\rho-\h{\psi}_0\{\tau(\bt),\bt\}]
    +[\h{\psi}_0\{\tau(\bt),\bt\}-\h{\psi}_0\{\tau(\hb),\hb\}]\\
    &=& O_p(n^{-1/2}).
  \end{eqnarray*}
  
  Using arguments similar to those in the proof of Proposition~\ref{decomp},
  \citet[Theorem~3.4]{huangb} extends to yield
  \[
  \h{\tau}(\hb)-\tau(\hb)
  = f_0\{\tau(\hb),\hb\}^{-1}\left[\rho-\h{\psi}_0\{\tau(\hb),\hb\}\right]
  +O_p(n^{-3/4}\log^{3/4}n),
  \]
  from which equation~(\ref{bth}) follows.

  By Taylor expansion,
  \begin{eqnarray*}
    \psi_0(\h{t},\hb)-\rho
    &=& \psi_0(\h{t},\hb)-\psi_0\{\tau(\hb),\hb\}\\
    &=&  f_0\{\tau(\hb),\hb\}\{\h{t}-\tau(\hb)\}+O_p(n^{-1}\log n),
  \end{eqnarray*}
  from which equation~(\ref{bsens}) follows. Equation~(\ref{bspec}) is
  established similarly.

\subsection*{Proof of Proposition~\ref{prop1}}
  The assertion follows by adapting the proof of \citet[corollary~2]{huangc}.

\subsection*{Proof of Corollary~\ref{coro1}}
  Let $\xbar{\phi}_\mathrm{cap}$ be the counterpart of $\h{\phi}_\mathrm{cap}$
  when $\xbar{\phi}(\cdot)$ is used in place of $\h{\phi}(\cdot)$.
  From Proposition~\ref{cra}, $\hb$ is also a near maximizer of
  $\xbar{\phi}({\bf b})$. Then, by adapting the proof of
  \citet[corollary~2]{huangc}, $n^{2/3}\{\xbar{\phi}_\mathrm{cap}-\phi(\hb)
  -\xbar{\phi}(\bt)+\phi(\bt)\}$ converges to a mean-zero, nondegenerate
  distribution. So is $n^{2/3}\{\h{\phi}_\mathrm{cap}-\phi(\hb)
  -\xbar{\phi}(\bt)+\phi(\bt)\}$ because $\h{\phi}_\mathrm{cap}=\xbar{\phi}_\mathrm{cap}+o_p(n^{-2/3})$ following \citet[theorem~3.4]{huangb}. The assertion
    follows
    \begin{eqnarray*}
      \psi_1(\h{\tau}_\mathrm{c},\hb) &=&
      \phi(\hb)-f_1\{\tau(\bt),\bt\}\{\h{\tau}_\mathrm{c}-\tau(\hb)\}
      +o_p(n^{-2/3})\\
      &=& \phi(\hb)-\lambda[\rho-\h{\psi}_0\{\tau(\bt),\bt\}]\\
      & & \mbox{}-f_1\{\tau(\bt),\bt\}\left(\h{\tau}_{\mathrm{c}}-\tau(\hb)-f_0\{\tau(\bt),\bt\}^{-1}
  [\rho-\h{\psi}_0\{\tau(\bt),\bt\}]\right)
  +o_p(n^{-2/3}),
    \end{eqnarray*}
    by equations~(\ref{bsens}) and~(\ref{bspec}), and Proposition~\ref{prop1}.

\subsection*{Proof of Corollary~\ref{coro41}}
  By extending \citet[lemma~3.2]{huangb}, it can be shown
  that there exists a constant $\epi>0$ such that
  \[
  \sup_{\|{\bf b}-\bt\|_\infty\leq\epi}|\h{\tau}({\bf b},\rho_n)
  -\tau({\bf b},\rho_n)|=
  O(n^{-1/2}\log^{1/2} n),
  \]
  almost surely. Since
  $\tau({\bf b},\rho_n)-\tau({\bf b})=O(n^{-1/2}\log^{1/2} n)$, it follows
  that
  \[
  \sup_{\|{\bf b}-\bt\|_\infty\leq\epi}|\h{\tau}({\bf b},\rho_n)-\tau({\bf b})|=
  O(n^{-1/2}\log^{1/2} n),
  \]
  almost surely. The proof of Theorem~\ref{thm2} then extends with only minor
  modifications.

\subsection*{Proof of Corollary~\ref{thmcip}}
  By the central limit theorem, $n_0^{1/2}[\h{\psi}_0\{\tau(\bt),\bt\}-\rho]$
  converges in distribution to a normal distribution. It follows that
  $n^{1/2}(\rho_n-\rho)$ converges to the $\delta$-quantile of this limiting
  distribution. Hence,
  $\rho_n-\rho=O(n^{-1/2})$ and equation~(\ref{bsens2}) applies. Meanwhile,
  in parallel with $\h{\tau}_\mathrm{c}=\h{\tau}(\hb)+O_p(n^{-2/3})$ following
  Theorem~\ref{thm2} and Proposition~\ref{prop1}, one can similarly show
  $\h{\tau}_\mathrm{c}(\rho_n)=\h{\tau}(\hb;\rho_n)+O_p(n^{-2/3})$.
  Consequently, equation~(\ref{bsens2}) remains valid with
  $\h{\tau}(\hb;\rho_n)$ replaced by $\h{\tau}_\mathrm{c}(\rho_n)$. Thus, the
  assertion follows.

\section*{Acknowledgements}
The author was supported in part by NIH Grants R01 CA230268,
R01 CA283687, and P30 AI050409.

\section*{Data availability statement}

The Breast Cancer Coimbra dataset analyzed in this article is openly available
at \url{https://archive.ics.uci.edu/dataset/451/breast+cancer+coimbra}
in UCI Machine Learning Repository, University of California, Irvine.

\end{document}